\documentclass[graybox]{svmult}
\usepackage{mathptmx,xcolor} 
\usepackage{helvet} 
\usepackage{courier} 
\usepackage{type1cm} 
\usepackage{makeidx} 
\usepackage{graphicx} 
\usepackage{multicol} 
\usepackage[bottom]{footmisc}
\usepackage{amsmath}
\usepackage{amssymb}
\usepackage{algorithm}
\usepackage{algpseudocode}
\makeindex 
\usepackage{caption}
\newtheorem{condition}{Condition}[section]
\begin{document}

\title*{Penalty Ensembles for Navier-Stokes with Random Initial Conditions and Forcing}
\titlerunning{Penalty Ensembles for Navier-Stokes}
\author{Rui Fang}

\institute{Rui Fang \at Department of Mathematics, University of Pittsburgh, Pittsburgh, PA 15260, \email{ruf10@pitt.edu}
}
\maketitle

\abstract*{In many applications, uncertainty in problem data leads to the need for numerous computationally expensive simulations. This report addresses this challenge by developing a penalty-based ensemble algorithm. Building upon Jiang and Layton's work on ensemble algorithms that use a shared coefficient matrix, this report introduces the combination of penalty methods to enhance its capabilities. Penalty methods uncouple velocity and pressure by relaxing the incompressibility condition. Eliminating the pressure results in a system that requires less memory. The reduction in memory allows for larger ensemble sizes, which give more information about the flow and can be used to extend the predictability horizon.}

\abstract{In many applications, uncertainty in problem data leads to the need for numerous computationally expensive simulations. This report addresses this challenge by developing a penalty-based ensemble algorithm. Building upon Jiang and Layton's work on ensemble algorithms that use a shared coefficient matrix, this report introduces the combination of penalty methods to enhance its capabilities. Penalty methods uncouple velocity and pressure by relaxing the incompressibility condition. Eliminating the pressure results in a system that requires less memory. The reduction in memory allows for larger ensemble sizes, which give more information about the flow and can be used to extend the predictability horizon.}

\section{Introduction}
\label{sec:1}
In fluid dynamics, uncertainties in problem data lead to the need for simulations with different initial and boundary conditions, body forces, and other data \cite{kalnay2003atmospheric}. They often result in computationally expensive simulations. The trajectories of these simulations can diverge significantly under slightly different initial conditions \cite{lorenz1963deterministic, lorenz1963predictability, lorenz1996essence}. The predictions are unreliable when they are $\mathbf{O}(1)$ apart.

We address these challenges using ensemble methods. Ensemble methods run multiple simulations with slightly different initial conditions. The best prediction is obtained by averaging the results from all simulations \cite{leith1974theoretical}.
Let $\Omega$ be an open regular domain in $\mathbb{R}^d$ ($d=2$ or 3) and $T$ be the final time. Consider the incompressible Navier--Stokes equations (NSE):
\begin{equation}\label{nse}
\begin{gathered}
\frac{\partial u}{ \partial t} + u\cdot \nabla u - \nu \Delta u +\nabla p = f, \text{ and } \nabla \cdot u =0 \text{ in } \Omega \times [0,T].
\end{gathered}
\end{equation}
Here, $u$ is the velocity, $p$ is the pressure, $f$ is the body force, and $\nu$ is the kinetic viscosity. $\nabla \cdot u$ denotes the divergence of the velocity field.

The velocity and pressure in equation (\ref{nse}) are coupled by the incompressible constraint $\nabla \cdot u =0$.
Penalty methods relax the incompressibility condition and uncouple velocity and pressure. One can eliminate the pressure through $p^\epsilon=-\frac{1}{\epsilon} \nabla \cdot u^\epsilon$ to speed up the calculation. The penalized NSE from \cite{temam1968methode} is described as the following:
\begin{equation}\label{penalty_conti}
\begin{gathered}
\frac{\partial u^\epsilon }{\partial t}+ u^\epsilon \cdot \nabla u^\epsilon + \frac{1}{2} (\nabla \cdot u^\epsilon) u^\epsilon -\nu \Delta u^\epsilon +\nabla p^{\epsilon} =f, \\
\quad \nabla \cdot u^{\epsilon}+\epsilon p^{\epsilon} =0, \text{ where } 0< \epsilon \ll 1.
\end{gathered}
\end{equation}
Temam \cite{temam1968methode} proved that $\lim_{\epsilon \to 0}(u^\epsilon, p^\epsilon)= (u, p)$.

Instead of a single realization, we conduct $J$ simulations. $u^\epsilon_j$ and $p^\epsilon_j$ represent the solutions to the NSE for the $j_{th}$ ensemble member and $f_j$ is the corresponding body force, where $1 \leq j \leq J$. We present the time discretization with a shared coefficient matrix for all ensemble members. We introduce the time step size $\Delta t$ and the final step $N= T/\Delta t$. We denote the velocity and pressure at the $n_{th}$ step as $u^{\epsilon,n}_j$ and $p^{\epsilon,n}_j$, respectively, where $0\leq n\leq N$.
\begin{definition}
Denote the ensemble mean and fluctuation at the $n_{th}$ step:
\begin{equation*}
\begin{gathered}
<u^\epsilon>^n: = \frac{1}{J} \sum_{j=1}^J u^{\epsilon,n}_j, \text{ and } U^{\epsilon,n}_j:= u^{\epsilon,n}_j - <u^\epsilon>^n.
\end{gathered}
\end{equation*}
\end{definition}
For clarity, we suppress the spatial discretization. We use an implicit-explicit time discretization and keep the resulting coefficient matrix independent of ensemble members. The method is to find $u^{\epsilon,n+1}_j\in X$, $p^{\epsilon,n+1}_j \in Q$:
\begin{equation}\label{BEFE-Penalty-Ensemble}
\begin{gathered}
\frac{u^{\epsilon, n+1}_j-u^{\epsilon, n}_j}{\Delta t} + <u^\epsilon>^n\cdot \nabla u^{\epsilon,n+1}_j + \frac{1}{2} (\nabla \cdot <u^\epsilon>^n) u^{\epsilon,n+1}_j \\
+U^{\epsilon, n}_j\cdot \nabla u^{\epsilon,n}_j + \frac{1}{2} (\nabla \cdot U^{\epsilon, n}_j) u^{\epsilon,n}_j
-\nu\Delta u_j^{\epsilon, n+1}
+\nabla p^{\epsilon, n+1}_{j} = f_j^{n+1},\\
\nabla\cdot u_j^{\epsilon, n+1} +\epsilon p_j^{\epsilon, n+1}= 0.
\end{gathered}
\end{equation} 
We can eliminate the pressure in equation (\ref{BEFE-Penalty-Ensemble}). Here $\epsilon$ is the same across all ensemble members to have a shared coefficient matrix. Thus, we assemble the matrix once and reduce the time of assembling the matrix separately for each ensemble member.
\section{Notation and preliminaries}\label{rec: preliminary}
The $L^2(\Omega)$ norm and the inner product are denoted by $\|\cdot\|$ and $(\cdot,\cdot)$. Likewise, the $L^p(\Omega)$ norms and the Sobolev $W^{k}_{p}(\Omega)$ norms are denoted by $\|\cdot\|_{L^p}$ and $\|\cdot\|_{W^{k}_{p}}$, respectively. $H^{k}(\Omega)$ is the Sobolev space $W^{k}_2(\Omega)$, with norm $\|\cdot\|_{k}$. For functions $v(t,x)$ defined for $t\in [0,T]$, we define, for $1\leq m< \infty$,
\begin{equation}
\|v\|_{\infty,k} := EssSup_{[0,T]} \|v(t,\cdot)\|_{k} \text{ and } \|v\|_{m,k} := \left(\int_{0}^T\|v(t,\cdot)\|^{m}_{k}\, dt \right)^{1/m}.
\end{equation}
Define the fluid velocity space and pressure space, respectively:
\begin{equation*}
\begin{gathered}
X:=(H^{1}_{0}(\Omega))^d=\{ v\in L^2(\Omega)^d:\nabla v\in L^2(\Omega)^{d\times d} \text{ and } v=0 \text{ on } \partial \Omega\}, \text{ and}\\
Q:=L^2_0(\Omega)=\{q\in L^2(\Omega): \int_{\Omega} q \, dx=0\}.
\end{gathered}
\end{equation*}
Let $X^h \subset X$ and $Q^h\subset Q$ be finite element spaces for fluid velocity and pressure respectively. We assume that $X^h$ and $Q^h$ satisfy the $LBB^h$ condition:
\begin{condition}\label{LBB_h}
Suppose $(X^h, Q^h)$ satisfies:
\begin{equation}
\inf_{ q^h \in Q^h} \sup_{v^h \in X^h} \frac{(q^h,\nabla \cdot v^h)}{\|v^h\|\|q^h\|}\geq \beta^h >0,
\end{equation}
where $\beta^h$ is bounded away from zero uniformly in $h$.
\end{condition}
We summarize standard properties of finite element spaces from \cite{john2016finite, layton2008introduction}. Assume $X^h$ and $Q^h$ satisfy the following approximation properties for $0\leq s\leq m$:
\begin{equation}\label{approximation-properties}
\begin{gathered}
\inf_{v\in X^h}\|u-v\|\leq Ch^{s+1} |u|_{s+1}, \forall u\in X \cap \left(H^{s+1}(\Omega)\right)^d,\\
\inf_{v\in X^h} \|\nabla (u-v)\|\leq Ch^s|u|_{s+1}, \forall u\in X \cap \left(H^{s+1}(\Omega)\right)^d,\\
\inf_{q\in Q^h}\|p-q\| \leq C h^s |p|_s, \forall p\in Q\cap H^s(\Omega).
\end{gathered}
\end{equation}
Denote the skew-symmetric trilinear form: $\forall u,v,w \in X$,
\begin{equation*}
\begin{split}
b^*(u,v,w) :=\frac{1}{2}(u\cdot \nabla v, w)-\frac{1}{2} (u\cdot \nabla w, v).
\end{split}
\end{equation*}
\begin{lemma}(\text{See \cite{layton2008introduction},  p.123 and p.155}) $\forall u,v,w \in X$, the trilinear term
$b^*(u,v,w)$ is equivalent to
\begin{equation*}
\begin{split}
b^*(u,v,w)=(u\cdot \nabla v, w)+\frac{1}{2} \left((\nabla \cdot u) v, w\right).
\end{split}
\end{equation*}
\end{lemma}
\begin{definition}\label{projection} $P_{Q^h}$ is the $L^2$ projection of $Q$ into $Q^h$. That is, for any $q\in Q$, $P_{Q^h}(q)$ satisfies
\begin{equation*}
(P_{Q^h} (q)-q , q^h)= 0
\ \forall q^h \in Q^h.
\end{equation*}
\end{definition}
\section{Algorithms}
Let $\Delta t_n$ denote the time step size at the $n_{th}$ step. The fully discrete approximation is given $(u^{\epsilon,n}_{j,h},p^{\epsilon,n}_{j,h}) \in (X^h, Q^h)$, find $(u^{\epsilon,n+1}_{j,h}, p^{\epsilon, n+1}_{j,h} ) \in (X^h, Q^h)$ satisfying:
\begin{equation}\label{method}
\begin{gathered}
\frac{(u^{\epsilon, n+1}_{j,h}-u^{\epsilon,n}_{j,h}, v^h)}{\Delta t_{n+1}} + b^*(<u^\epsilon_h>^n, u^{\epsilon,n+1}_{j,h},v^h)+ b^*(u^{\epsilon,n}_{j,h} -<u^\epsilon_h>^n, u^{\epsilon,n}_{j,h},v^h)\\+
\nu (\nabla u^{\epsilon,n+1}_{j,h}, \nabla v^h)-(p^{\epsilon,n+1}_{j,h}, \nabla \cdot v^h) = (f_j^{n+1},v^h),\\
(\nabla \cdot u^{\epsilon, n+1}_{j,h}, q^h)
+ \epsilon (p^{\epsilon, n+1}_{j,h}, q^h)=0.
\end{gathered}
\end{equation}
for all $(v^h, q^h)\in (X^h,Q^h)$. We reduce computation time by solving the ensemble members with a shared coefficient matrix and different right-hand side vectors. We can eliminate the pressure. The system will require less memory. The penalty-based ensemble method is presented in Algorithm \ref{alg:cap}.
\begin{algorithm}
\caption{Penalty-based ensembles}\label{alg:cap}
\begin{algorithmic}
\Require $\Delta t_{n+1}$, $\epsilon$, $u^{\epsilon,n}_{j,h}, j=1,\ldots, J$ at time $t$.
\While{$t < T$}:
\State Compute the average $<u^\epsilon_h>^n = \frac{1}{J} \sum_{j=1}^J u^{\epsilon,n}_{j,h}$.
\State Solve the next step velocity $u^{\epsilon, n+1}_{j,h}$ by the method in equation (\ref{method}).
\State Update the average $<u^\epsilon_h>^{n+1} = \frac{1}{J} \sum_{j=1}^J u^{\epsilon,n+1}_{j, h}$.
\State Calculate $\max_{j} \|\nabla (u^{\epsilon, n+1}_{j,h} - <u^\epsilon_h>^{n+1})\|$ and verify the CFL condition.
\If{the CFL condition is satisfied}
\State update time by $t = t + \Delta t_{n+1}$.
\Comment{Succeed, and go to next timestep}
\Else
\State $\Delta t_{n+1} = \frac{1}{2} \Delta t_{n+1}$.
\EndIf
\EndWhile
\end{algorithmic}
\end{algorithm}
\section{Error Analysis}
We summarize stability and convergence theorems from our analysis, Fang \cite{fang2024numerical_analysis}. Due to the explicit discretization of the stretching term
$b^*(u^{\epsilon,n}_{j,h} -<u^\epsilon_h>^n, u^{\epsilon,n}_{j,h},v^h)$, a CFL condition \cite{courant1928partiellen} is necessary to ensure its stability and convergence of the method described by equation (\ref{method}).
\begin{theorem}\label{thm: stability} (Fang \cite{fang2024numerical_analysis}) Consider the method in equation (\ref{method}). Suppose for any $\ 0 \leq n\leq N-1$,
\begin{equation*}
C\frac{\Delta t} {\nu h} \|\nabla (u^{\epsilon,n+1}_{j,h}- <u^\epsilon_h>^{n+1})\|^2 \leq 1.
\end{equation*}
Then for any $N>1$, the following holds
\begin{align*}
\begin{gathered}
\frac{1}{2} \|u^{\epsilon, N}_{j,h}\|^2 + \frac{1}{4}\sum_{n=0}^{N-1} \|u^{\epsilon,n+1}_{j,h}- u^{\epsilon,n}_{j,h}\|^2 + \frac{\nu \Delta t}{4} \|\nabla u^{\epsilon,N}_{j,h}\|^2\\+ \frac{\nu\Delta t}{4}\sum_{n=0}^{N-1}\|\nabla u^{\epsilon,n+1}_{j,h}\|^2
+ \frac{\Delta t}{\epsilon}\sum_{n=0}^{N-1} \|P_{Q^h}(\nabla \cdot u^{\epsilon,n+1}_{j,h})\|^2
\\
\leq \frac{\Delta t}{2 \nu} \sum_{n=0}^{N-1} \|f^{n+1}_{j}\|^2_{-1} +\frac{1}{2} \|u^\epsilon_{j,h}(0)\|^2 + \frac{\nu \Delta t}{4}\|\nabla u^\epsilon_{j,h}(0)\|^2.
\end{gathered}
\end{align*}
\end{theorem}
The error of $j_{th}$ ensemble member at time $t_n$ is denoted by
\begin{equation*}e^{\epsilon,n}_{j,h} := u^{\epsilon}_j(t_n) - u^{\epsilon,n}_{j,h}.
\end{equation*}
\begin{theorem}(Fang \cite{fang2024numerical_analysis})\label{convergence_ensemble_penalty}
Let $\Omega$ be a convex polygonal/polyhedral domain. Consider the method in equation (\ref{method}). Under the assumption of Theorem \ref{thm: stability}, there exists  positive constants $C$ and $C_0$ independent of $h$ and $\Delta t$ such that 
\begin{equation}\label{convergence_eqn}
\begin{gathered}
\|e^{\epsilon, N}_{j,h}\|^2 + \frac{1}{2}\sum_{n=0}^{N-1} \|e^{\epsilon,n+1}_{j,h}- e^{\epsilon,n}_{j,h}\|^2 + \Delta t \nu \|\nabla e^{\epsilon,N}_{j,h}\|^2\\
+ C_0\Delta t \sum_{n=0}^{N-1}\nu\|\nabla e^{\epsilon,n+1}_{j,h}\|^2 \leq \exp (\alpha )\Biggl\{\|e^{\epsilon, 0}_{j,h}\|^2+ \Delta t \nu \|\nabla e^{\epsilon,0}_{j,h}\|^2\\
+ h^{2m}C(\nu) T\left( \||u^{\epsilon}_{j,t}|\|^2_{\infty,0} + \frac{1}{\nu^2} \||p^{\epsilon}_{j,t}|\|^2_{\infty,0} \right)+ (\Delta t)^3C(\nu) \||u^{\epsilon}_{j,tt}|\|^2_{\infty, 0}\\
+ h^{2m}\epsilon \Delta t C(\nu, \beta^h)\left(\||u^{\epsilon}_{j,t}|\|^2_{2,0}  + \||p^{\epsilon}_{j,t}|\|^2_{2,0} \right)\\
+ h^{2m}C(\nu) T\left(  \|| u^{\epsilon}_j|\|^2_{2,0} +\frac{1}{\nu^2} \||p^{\epsilon}_j|\|^2_{2,0}\right)+ C(\nu) (\Delta t)^2 \||\nabla u^{\epsilon}_{j,t}|\|^2_{\infty,0}
\Biggl\},
\end{gathered}
\end{equation}
where
\begin{equation*}
\alpha = C(\nu) \Delta t \sum_{n=0}^{N-1} \|\nabla u^{\epsilon, n+1}\|^4. 
\end{equation*}
\end{theorem}
Combining Theorem \ref{convergence_ensemble_penalty} with the result of Shen \cite{shen1995error}, Theorem $4.1$, p. 395, we have the following.
\begin{corollary} (Fang \cite{fang2024numerical_analysis})
Under the assumption of Theorem \ref{thm: stability}, for regular solutions, we have the following optimal estimates:
\begin{equation*}
\max_{t_n} \| u_j(t_n) - u^{\epsilon,n}_{j,h}\|^2 +\Delta t \sum_{n=1}^{N} \|\nabla ( u_j(t_n) - u^{\epsilon,n}_{j,h})\|^2 \leq C(u,\nu) (\epsilon +\Delta t + h^m)^2.
\end{equation*}
\end{corollary}

\section{Numerical Experiments}
\label{sec:2}
We conduct numerical experiments with simply $J=2$ ensemble members. In the first test, we applied perturbations on the initial conditions of the modified Green--Taylor vortex. We verified the accuracy of the algorithm and confirmed the predicted convergence rates. We studied a rotating flow on offset cylinders and adapted the time step in the second test. The kinetic energy and enstrophy suggested that our algorithm preserved the stability.

We employ the Taylor--Hood ($P2-P1$) finite element pair. Here, $P2$ represents a second-order polynomial for the velocity, while $P1$ represents a first-order polynomial for the pressure. The unstructured mesh was generated with GMSH, an open-source finite element mesh generator \cite{geuzaine2009gmsh}.
\subsection{Convergence Experiment}
\label{subsec:2}
The Green--Taylor vortex problem is commonly used for 
calculating convergence rates, since the analytical solution is known \cite{john2016finite}. In $\Omega = (0,1)^2$, the exact solution of the modified Green--Taylor vortex is
\begin{align*}
\begin{gathered}
u(x,y,t) = (-\cos(x) \sin(y) \sin(t), \sin(x) \cos(y) \sin(t))^\top,\\
p(x,y,t) = \frac{1}{4}\left( \cos(2x)+ \cos(2y)\right) \sin^2(t).
\end{gathered}
\end{align*}
We take the exact initial condition $u(x,y,0)$ and perturb for the ensemble members with $\delta_1 = 10^{-3}$ and $\delta_2 = -10^{-3}$:
\begin{align*}
\begin{gathered}
u_1(x,y,0) = (1+\delta_1) u(x,y,0),\ u_2(x,y,0) = (1+ \delta_2) u(x,y,0).
\end{gathered}
\end{align*}
We set the kinematic viscosity $\nu = 1$, the characteristic velocity of the flow $U=1$, the characteristic length $L=1$, and the Reynolds number $Re = \frac{UL}{\nu}$. To discretize the domain, we choose a sequence of mesh sizes $h=\frac{1}{g}$, where $g$ is successively defined as $g = \left(\frac{3}{2}\right)^i \cdot 27$, where $i=0, \ldots, 4$. We set the time step size $\Delta t = \frac{h}{10}$, penalty parameter $\epsilon= \Delta t$, and the final time $T=1$. Denote the error $e(h) = C h^\beta$, solve $\beta$ via
\begin{equation*}
\beta = \frac{\ln(e(h_1)/e(h_2))}{\ln(h_1/h_2)}
\end{equation*}
at two successive values of $h$, where $\beta$ is the convergence rate. Set inhomogeneous Dirichlet boundary condition $
u_h = u_{true} \text{ on } \partial \Omega$ for all ensemble members.

\begin{table}
\caption{Errors and convergence rates for $u_1$.}
\label{tab:1}
\centering
\begin{tabular}{c | c | c | c | c }
g & $\max_{t_n} \| u_1(t_n) - u^{\epsilon,n}_{1,h}\|$ & rate & $\sqrt{\Delta t \sum_{n=1}^{N} \|\nabla ( u_1(t_n) - u^{\epsilon,n}_{1,h})\|^2} $ & rate \\
\hline
$(\frac{3}{2})^0 \cdot27$ & $1.38\cdot 10^{-4}$ & -- & $3.61\cdot 10^{-4}$ & -- \\
$(\frac{3}{2})^1\cdot 27$ & $9.37\cdot 10^{-5}$ &0.99 & $2.38\cdot 10^{-4}$&1.05 \\
$(\frac{3}{2})^2 \cdot 27$ & $6.26\cdot 10^{-5}$ &0.99& $1.57\cdot 10^{-4}$ &1.02\\
$(\frac{3}{2})^3 \cdot 27$ & $4.14\cdot 10^{-5}$ &0.99 & $1.03\cdot 10^{-4}$ & 1.01\\
$(\frac{3}{2})^4 \cdot 27$ & $2.78 \cdot 10^{-5}$ & 0.99 & $6.90\cdot 10^{-5}$ & 1.00
\end{tabular}
\end{table}

\begin{table}
\caption{Errors and convergence rates for $u_2$.}
\label{tab:2}
\centering
\begin{tabular}{c | c | c | c | c }
g& $\max_{t_n} \| u_2(t_n) - u^{\epsilon,n}_{2,h}\|$ & rate & $\sqrt{\Delta t \sum_{n=1}^{N} \|\nabla ( u_2(t_n) - u^{\epsilon,n}_{2,h})\|^2}$ & rate \\
\hline
$(\frac{3}{2})^0 \cdot27$ & $1.38\cdot 10^{-4}$ & -- & $3.59\cdot 10^{-4}$ &-- \\
$(\frac{3}{2})^1\cdot 27$ & $9.34\cdot 10^{-5}$ &0.99 & $2.38\cdot 10^{-4}$ & 1.06\\
$(\frac{3}{2})^2 \cdot 27$ & $6.24\cdot 10^{-5}$ & 0.99& $1.56 \cdot 10^{-4}$ & 1.02\\
$(\frac{3}{2})^3 \cdot 27$ & $4.13 \cdot 10^{-5}$ &0.99& $1.02\cdot 10^{-4}$ & 1.01 \\
$(\frac{3}{2})^4 \cdot 27$ & $2.76\cdot 10^{-5}$ & 0.99 &$6.87\cdot 10^{-5}$ & 1.00
\end{tabular}
\end{table}
Table~\ref{tab:1} and Table~\ref{tab:2} show that the convergence rates of $u_1$ and $u_2$ are first order as predicted.

\subsection{Stability Verification}
\label{subsec:3}
The domain is a disk with a smaller off-center obstacle inside. Let outer circle radius $r_1 = 1$ and inner circle radius $r_2 = 0.5$, $c = (c_1, c_2)= (\frac{1}{2},0)$. The domain is
\begin{equation*}
\Omega = \{ (x,y): x^2 + y^2 \leq r_1^2 \text{ and } (x-c_1)^2 + (y-c_2)^2 \geq r_2^2\}.
\end{equation*}
A counterclockwise rotational force drives the flow,
\begin{equation*}
f(x,y,t) = (-4 y(1-x^2-y^2), 4x (1-x^2-y^2))^\top,
\end{equation*}
with no-slip boundary conditions. Note that the outer circle remains
stationary. We set mesh size $h = 0.05$. We chose
the final time $T = 100$, time step size $\Delta t = h/10$, penalty parameter $\epsilon = \Delta t$, $\nu = 1/150, L=1, U=1$, and $Re = \frac{UL}{\nu}$.

Initial condition $u(x,y,0) = 0$ and the Dirichlet boundary condition $u=0$ on $\partial \Omega$. We denote two perturbation values $\delta_1 = 0.1$ and $\delta_2 = -0.1$ to introduce perturbations. The perturbed initial conditions for $u_1$ and $u_2$:
\begin{equation*}
u_j(x,y,0) := \delta_j (1-x^2-y^2)(0.5^2 - (x-0.5)^2 -y^2), \text{where } j=1, 2.
\end{equation*}
This formulation preserves the no-slip boundary condition. For the evaluation, we calculated flow statistics, enstrophy, and kinetic energy:
\begin{equation*}
\begin{gathered}
\text{enstrophy} = \frac{1}{2} \nu \|\nabla \times \bar{u}\|^2 \text{ and }
\text{kinetic energy} = \frac{1}{2} \| \bar{u} \|^2.
\end{gathered}
\end{equation*}
To ensure the algorithm is stable, we adopt $C = 1200$ from \cite{jiang2014algorithm} and adapt $\Delta t$ to satisfy 
\begin{equation}\label{CFL}
\frac{\Delta t}{ h} \|\nabla (u^{\epsilon,n+1}_{j,h}- <u^\epsilon_h>^{n+1})\|^2 \leq \frac{1200}{Re}, \text{ where } j=1,2.
\end{equation}
If equation (\ref{CFL}) is violated, we reduce $\Delta t $ by half and repeat the process until equation (\ref{CFL}) is satisfied.
\begin{figure}[h!]
\centering
\includegraphics[width=\textwidth]{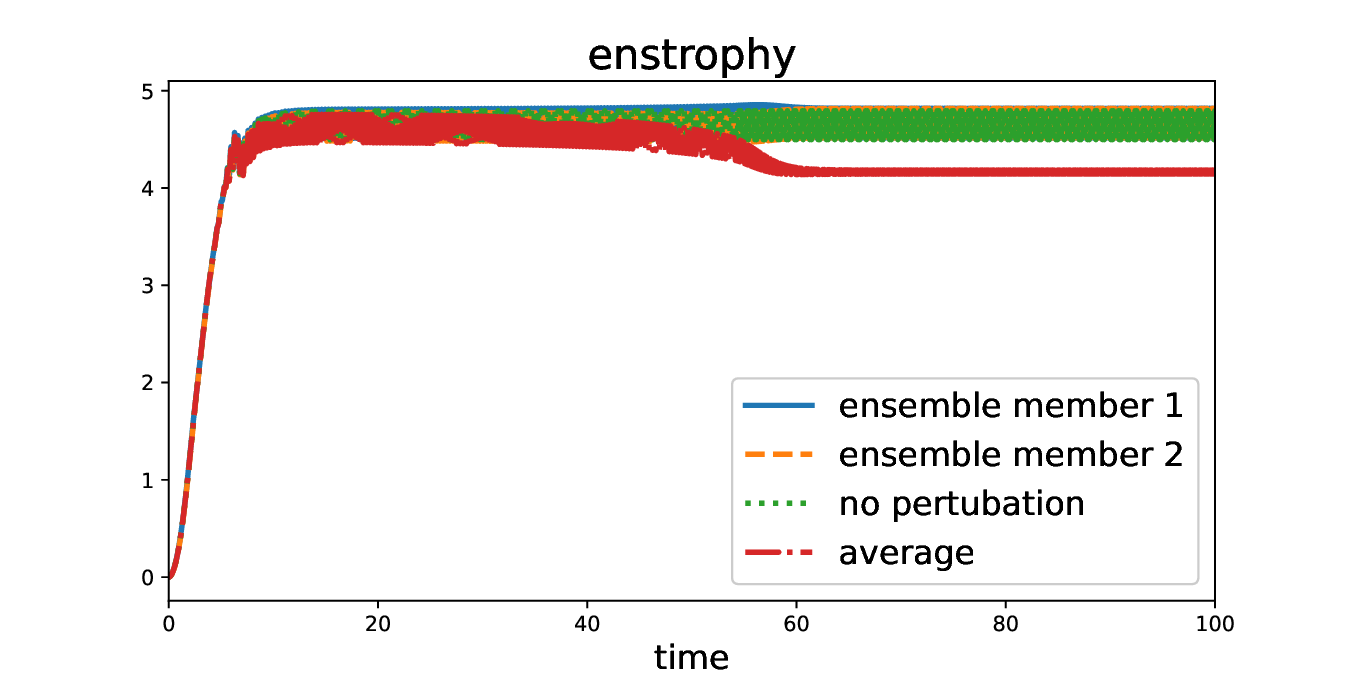}
\caption{The enstrophy of the ensemble average $u_{ave}$ shows that the flow is smoothed out on average.}
\label{fig: enstrophy}
\end{figure}
\begin{figure}[h!]
\centering
\includegraphics[width=\textwidth]{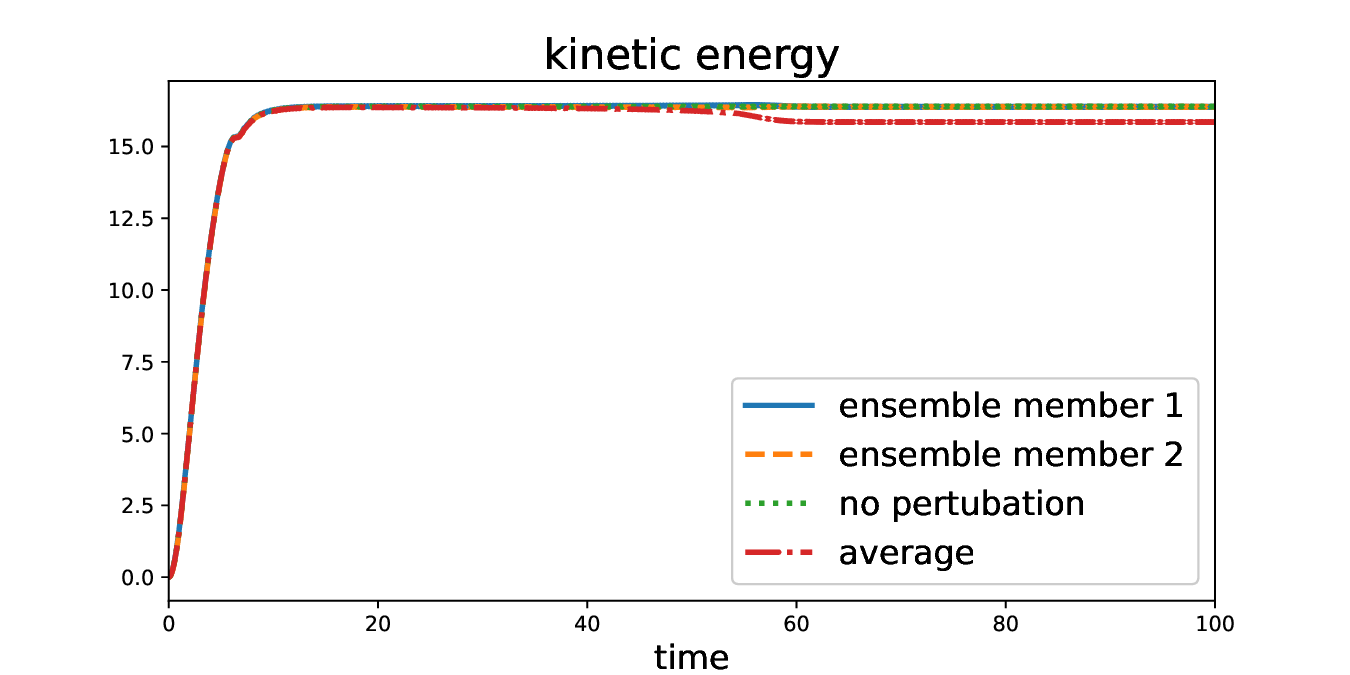}
\caption{The kinetic energy is similar across all members and has reached a statistically steady state.}
\label{fig:kinetic-energy}
\end{figure}

In Figure \ref{fig: enstrophy}, the enstrophy of ensemble average $u_{ave}$ is lower than the ensemble members and no perturbed one, which shows that the flow is smoothed out on average. In Figure \ref{fig:kinetic-energy}, the kinetic energy is similar across all members and has reached a statistically steady state.
\begin{figure}[h!]
\centering
\includegraphics[width =\linewidth]{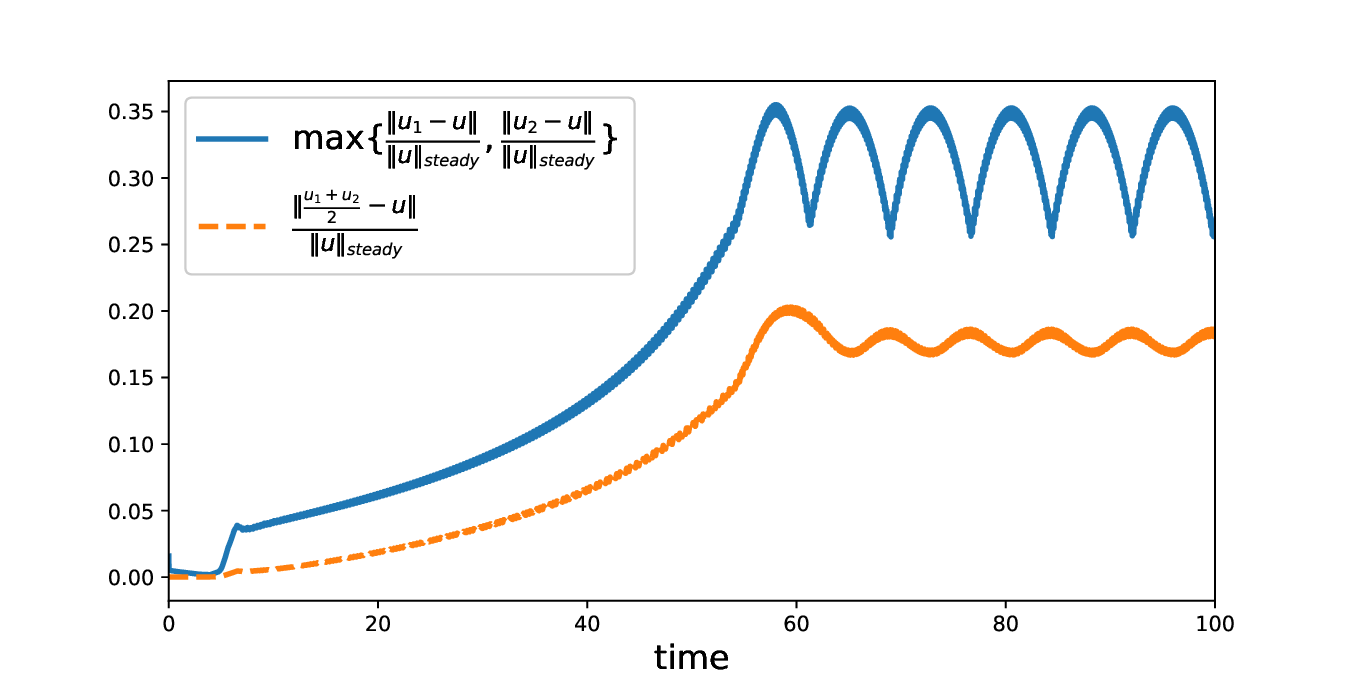}
\caption{The plot of the relative error $\|u_j - u\|/\|u\|_{steady}$.}
\label{fig: spread}
\end{figure}

Figure \ref{fig: spread} plots the relative error $\|u_j - u\|/\|u\|_{steady}$, the horizontal axis is the time, and the vertical axis is the relative error. We normalized by $\|u\|$ after the norm of the flow velocity reaches a statistically steady state. The predictability horizon is defined as the duration until $\|u_j - u\|/\|u\|_{steady}$ reaches a chosen threshold $0.1$. Computing the maximum error at each time step for ensemble members 1 and 2, we found that the $0.1$-predictability suggests that the predictability horizon for a single realization is $T=32.6$, while $T=47.5$ for the ensemble average. When we set the threshold to approximately $0.2$, it suggests that the ensemble average remains accurate until the final time $T=100$, which is significantly longer than a single realization. Our algorithm uses a shared coefficient matrix, reducing assembling time. Eliminating pressure reduces memory. It allows for a larger size of ensemble members. Large ensemble sizes improve statistical robustness, reduce sampling error, better represent uncertainty and thus yield more reliable results.

\section{Conclusions}
In this report, we test the stability and convergence of the penalty-based ensemble method in equation (\ref{method}). Our algorithm uses a shared coefficient matrix and reduces the memory, which allows for larger ensemble sizes. It gives a more accurate and reliable simulation. The ensemble spread is related to the forecast uncertainty, but the nature of this relationship remains unclear and is ongoing research \cite{van2016probabilistic}. Open problems include extending this to a higher Reynolds number and adapting penalty parameters \cite{fang2024numerical}.
\begin{acknowledgement}
I thank Professor William Layton for his instruction and support throughout this research. This research herein of Rui Fang was supported in part by
the NSF under grant DMS 2110379 and DMS 2410893.
\end{acknowledgement}
\bibliographystyle{spmpsci}
\bibliography{mybib.bib}
\end{document}